\documentstyle[amscd,amssymb,verbatim,12pt]{amsart}
\newcommand{\del}{\partial}
\newcommand{\ssl}{\frak{sl}}
\newcommand{\ch}{\operatorname{ch}}
\newcommand{\CH}{\operatorname{CH}}

\newcommand{\OO}{{\cal O}}

\newcommand{\DD}{{\cal D}}

\newcommand{\SL}{\operatorname{SL}}

\newcommand{\Th}{\Theta}

\newcommand{\Pic}{\operatorname{Pic}}

\newcommand{\de}{\delta}

\numberwithin{equation}{section}

\newtheorem{thm}{Theorem}[section]
\newtheorem{prop}[thm]{Proposition}
\newtheorem{lem}[thm]{Lemma}
\newtheorem{cor}[thm]{Corollary}
\newenvironment{rem}{\vspace{3mm}\noindent
{\bf Remark.}}{\vspace{3mm}}

\newenvironment{rems}{\vspace{3mm}
\noindent {\bf Remarks.}}{\vspace{3mm}}

\newcommand{\Pf}{\noindent {\it Proof}}
\newcommand{\id}{\operatorname{id}}

\newcommand{\ra}{\rightarrow}

\newcommand{\TT}{{\cal T}}

\newcommand{\PP}{{\cal P}}

\newcommand{\LL}{{\cal L}}

\newcommand{\End}{\operatorname{End}}

\newcommand{\De}{\Delta}
\newcommand{\la}{\lambda}
\newcommand{\th}{\theta}

\newcommand{\Z}{{\Bbb Z}}
\newcommand{\Q}{{\Bbb Q}}

\newcommand{\wt}{\widetilde}

\newcommand{\ad}{\operatorname{ad}}
\newcommand{\sub}{\subset}
\newcommand{\ed}{\qed\vspace{3mm}}

\title{Lie symmetries of the Chow group of a Jacobian and the tautological subring}
\author{A. Polishchuk}
\address{Department of Mathematics, University of Oregon, Eugene, OR 97405}
\email{apolish@@math.uoregon.edu}
\thanks{Supported in part by NSF grant DMS-0302215}

\begin{document}
\begin{abstract} Let $J$ be the Jacobian of a smooth projective curve.
We define a natural action of the Lie algebra of polynomial Hamiltonian vector
fields on the plane, vanishing at the origin, on the Chow group $\CH(J)_{\Q}$.
Using this action we obtain some relations between tautological cycles in $\CH(J)_{\Q}$.
\end{abstract}
\maketitle

\bigskip

\centerline{\sc Introduction}

\medskip

Let $J$ be the Jacobian of a smooth projective curve $C$ of genus $g\ge 2$.
We fix a point $x_0\in C$ and consider the corresponding embedding
$$\iota:C\ra J$$
mapping a point $x\in C$ to the isomorphism class of the line bundle
$\OO_C(x-x_0)$. We always consider $C$ as a subvariety in $J$ via
this embedding.
We define the {\it tautological subring} 
$\TT\CH(J)_{\Q}\sub\CH(J)_{\Q}$ in the Chow ring of $J$ with coefficients in $\Q$
as the smallest subring containing the class $[C]$ of the curve and closed
under taking pull-backs with respect to the
natural isogenies $[n]:J\to J$ and under the Fourier transform 
$S:\CH(J)_{\Q}\to\CH(J)_{\Q}$ (defined in \cite{B1}).
The corresponding subring in the quotient of $\CH(J)_{\Q}$ modulo algebraic
equivalence was considered by Beauville in \cite{Bmain} and by the author in \cite{P}.
It is known that modulo algebraic equivalence this subring is generated by the characteristic
classes of the Picard bundle on $J$. Also, a number of nontrivial relations between
these generators (still modulo algebraic equivalence) was described in \cite{P}.
In the present paper we will show how to lift these relations to the Chow ring
(after adding some more generators). 
This is achieved using the action of a certain Lie algebra on
$\CH(J)_{\Q}$ 
extending the well action of $\ssl_2$ associated with the natural polarization of $J$
(see \cite{K}). The construction of this action
may be of independent interest.

To state the results precisely we need to introduce some more notation.
Recall that the Chow ring of $J$ with rational coefficients admits a decomposition
$$\CH(J)_{\Q}=\oplus_{p,s}\CH^p_s(J),$$
where $\CH^p_s$ consists of $c\in\CH^p(J)_{\Q}$ such that $[n]^*c=n^{2p-s}c$
(see \cite{B2}). For every class $c\in\CH^p(J)_{\Q}$ we denote by $c_s\in\CH^p_s(J)$ 
its components with respect to the above decomposition. We also denote by
$\th\in\CH^1(J)_{\Q}$ the class of a symmetric theta divisor.

Let us define two families of classes in $\CH(J)_{\Q}$ by setting
$$p_n=S([C]_{n-1}),\  n\ge 1,$$
$$q_n=S(\th\cdot [C]_n),\ n\ge 0.$$
Note that $p_1=-\th$ and $q_0=g\cdot [J]$. It is easy to see that all the classes $(p_n)$ and $(q_n)$
belong to the tautological subring $\TT\CH(J)_{\Q}$. All the classes $q_n$ for $n\ge 1$ are
algebraically equivalent to zero.

\begin{thm}\label{com-rel-thm} 
There exists a family of operators $(X_{m,n}, Y_{m,n})$ on $\CH(J)_{\Q}$,
where $m,n\in \Z$, such that $X_{m,n}=0$ unless $m,n\ge 0$ and $m+n\ge 2$ (resp.,
$Y_{m,n}=0$ unless $m,n\ge 0$), satisfying the 
commutation relations
\begin{equation}\label{main-com-rel}
[X_{m,n},X_{m',n'}]=(nm'-mn')X_{m+m'-1,n+n'-1},
\end{equation}
\begin{equation}\label{main-com-rel2}
[X_{m,n},Y_{m',n'}]=(nm'-mn')Y_{m+m'-1,n+n'-1},
\end{equation}
$$[Y_{m,n},Y_{m',n'}]=0$$
and such that
$$\frac{1}{n!}X_{0,n}(a)=p_{n-1}\cdot a,$$
$$\frac{1}{n!}Y_{0,n}(b)=q_n\cdot a.$$
Furthermore, one has
\begin{equation}\label{FourXmn}
SX_{m,n}S^{-1}=(-1)^n X_{n,m},
\end{equation}
\begin{equation}\label{FourYmn}
SY_{m,n}S^{-1}=(-1)^n Y_{n,m}.
\end{equation}
\end{thm}

\begin{rem} In fact, one can easily see from the proof that the above operators
on the Chow group are induced by endomorphisms of a $\Q$-motive of $J$ and the relations
are satisfied already for this motive action. 
\end{rem}

Explicit formulas for operators $X_{m,n}$ and $Y_{m,n}$ will be given in section \ref{com-rel-sec}
(see \eqref{Y-eq} and \eqref{X-eq}).
We will also show that the tautological subring $\TT\CH(J)_{\Q}$ is closed under all operators
$X_{m,n}$ and $Y_{m,n}$.
Note that the commutation relation \eqref{main-com-rel} is the defining relation for
the Lie algebra of polynomial Hamiltonian vector fields on the plane with the standard
symplectic form (see e.g., \cite{Fuks}, ch.~1,\S 1). The restriction $m+n\ge 2$ that we imposed for our generators $X_{m,n}$ corresponds to considering the subalgebra of vector fields vanishing
at the origin. Also, note that 
the operators $(X_{2,0}/2,X_{1,1},X_{0,2}/2)$ generate the well known action of $\ssl_2$
on $\CH(J)_{\Q}$ (see \cite{K}).

As an application of the above Lie action we prove the following result.

\begin{thm}\label{mainthm} 
(i) The ring $\TT\CH(J)_{\Q}$ is generated by the classes $(p_n)$ and $(q_n)$.
Furthermore, let us consider the following differential operator 
$$\DD=\frac{1}{2}\sum_{m,n\ge 1}{m+n\choose n}p_{m+n-1}\del_{p_m}\del_{p_n}+
\sum_{m\ge 1,n\ge 1}{m+n-1\choose n}q_{m+n-1}\del_{q_m}\del_{p_n}-
\sum_{n\ge 1}q_{n-1}\del_{p_n},$$
where $(\del_{p_n})$ (resp., $(\del_{q_m})$) are partial derivatives with respect to $(p_n)$
(resp., $(q_n)$).
Then the space of polynomial relations between $(p_n,q_n)$ in $\CH(J)_{\Q}$ is stable
under the action of $\DD$.

(ii) The operators $X_{m,n}$ and $Y_{m,n}$ preserve the subspace $\TT\CH(J)_{\Q}$ and act on it via the following differential operators (for $m\ge 1$):
\begin{align*}
&\frac{(-1)^m}{m!}X_{m,n}|_{\TT\CH(J)_{\Q}}=
\frac{1}{m!}\sum_{i_1,\ldots,i_m\ge 1}\frac{(n+i_1+\ldots+i_m)!}{i_1!\ldots i_m!}p_{n+i_1+\ldots+i_m-1}
\del_{p_{i_1}}\ldots\del_{p_{i_m}}\\
&+\frac{1}{(m-1)!}\sum_{i_1,\ldots,i_{m-1};j\ge 1}\frac{(n+i_1+\ldots+i_{m-1}+j-1)!}{i_1!\ldots 
i_{m-1}! (j-1)!}q_{n+i_1+\ldots+i_{m-1}+j-1}
\del_{p_{i_1}}\ldots\del_{p_{i_{m-1}}}\del_{q_j}\\
&-\frac{1}{(m-1)!}\sum_{i_1,\ldots,i_{m-1}\ge 1}\frac{(n+i_1+\ldots+i_{m-1})!}{i_1!\ldots 
i_{m-1}!}q_{n+i_1+\ldots+i_{m-1}-1}
\del_{p_{i_1}}\ldots\del_{p_{i_{m-1}}},
\end{align*}
where for $m=1$ the last term should be understood as $-n!q_{n-1}$, 
$$(-1)^mY_{m,n}|_{\TT\CH(J)_{\Q}}=\sum_{i_1,\ldots,i_m\ge1}
\frac{(n+i_1+\ldots+i_m)!}{i_1!\ldots i_m!}q_{n+i_1+\ldots+i_m}\del_{p_{i_1}}\ldots\del_{p_{i_m}}.$$
\end{thm}

It is easy to check that the above differential operators in independent variables $(p_n,q_n)$
satisfy relations \eqref{main-com-rel} and \eqref{main-com-rel2}. Since $\DD$ corresponds to the action of $X_{2,0}/2$, it follows that any $\DD$-invariant ideal in $\Q[p,q]$
is also invariant under all the other differential operators above.
Since the operator $\DD$ lowers the degree by $1$ (where $\deg p_n=\deg q_n=n$),
starting from vanishing of polynomials of degree $g+1$ 
and applying powers of $\DD$ we get nontrivial relations
in $\CH(J)_{\Q}$ (see section \ref{ex-sec} for some examples).
Since all the classes $q_n$ for $n\ge 1$ are algebraically equivalent to zero, 
we recover the relations between $(p_n)$ modulo algebraic equivalence proved in \cite{P}.
The above theorem also has the following corollary closely related to the work of Beauville
\cite{B0}. Consider the group of $0$-cycles $\CH^g(J)_{\Q}$ equipped with the Pontryagin product.
Let $K$ be the canonical class on the curve $C$. Then we have a special $0$-cycle
$\iota_*K\in\CH^g(J)_{\Q}$. The proof of the following corollary will be given in section \ref{prelim-sec}.

\begin{cor}\label{zero-cor}
The intersection $\TT\CH(J)_{\Q}\cap\CH^g(J)_{\Q}$ coincides with the $\Q$-subalgebra with
respect to the Pontryagin product generated by the classes $[n]_*\iota_*K$, where $n\in\Z$.
\end{cor}

\begin{rems} 1. It is easy that the classes of the subvarieties of special divisors $W_d\sub J$ 
belong to $\TT\CH(J)_{\Q}$. More precisely, we will show in section \ref{ex-sec} that they can
be expressed as universal polynomials in classes $(p_n-q_n)$.

\noindent 2. Of course, the ring $\TT\CH(J)_{\Q}$ depends on a choice of a point $x_0\in C$.
For example, if $(2g-2)x_0=K$ in $\Pic(C)_{\Q}$ then all $q_n$ vanish. In fact,
using Abel's theorem it is easy to see that $(2g-2)x_0=K$ iff $q_1=0$. The vanishing of other
classes $q_n$ in this case follows also from the formula
$$q_n=\DD(q_1p_n)+q_1q_{n-1}.$$
 \end{rems}


\noindent{\it Notation}. 
We use the convention ${n\choose m}=0$ for 
$m<0$ and for $n<m$.

\section{Preliminaries}\label{prelim-sec}

Let $\Th\sub J$ by a symmetric theta divisor (corresponding to some choice of a
theta characteristic on $C$), so that $\th=c_1(\Th)$. 
Consider the line bundle on $J\times J$ given by
$\LL=\OO_{J\times J}(p_1^{-1}\Th+p_2^{-1}\Th-m^{-1}\Th)$
where $p_1,p_2:J\times J\to J$ are the natural projections and
$m=p_1+p_2:J\times J\to J$ is the group law.
It is easy to see that $\LL|_{C\times C}\simeq\OO_C(\De_C-x_0\times C-C\times x_0)$,
where $\De_C\sub C\times C$ is the diagonal. Indeed, it suffices to check
that $\LL|_{C\times J}\simeq\PP_C$, where $\PP_C$ is the universal family of degree $0$
line bundles on $C$ trivialized at $x_0$. This in turn follows from the fact that $\LL^{-1}$ corresponds to
the normalized Poincar\'e line bundle on $J\times \hat{J}$ under the principal polarization
isomorphism $\phi:J\wt{\ra}\hat{J}$ and from the equality $\phi\circ \iota=-a$, where
$a:C\to\hat{J}$ is the embedding induced by $\PP_C$ (see \cite{P-av}, 17.3).
We denote by $S$ the Fourier transform on $\CH(J)_{\Q}$ defined by
$$S(c)=p_{2*}(\exp(c_1(\LL))\cdot p_1^*c).$$
We refer to \cite{B1} for the detailed study of this transform.
In particular we will use the following properties:
$$S^2=(-1)^g[-1]^*,$$
$$S(\CH^p_s(J))\sub\CH^{g-p+s}_s(J),$$
$$S(a*b)=S(a)\cdot S(b),$$
where $a*b$ denotes the Pontryagin product on $\CH(J)_{\Q}$.

It is easy to see that $S([C])=\sum_{n\ge 1}p_n$ is exactly the decomposition of the class
$S([C])$ into components of different codimensions, so that $p_n\in\CH^n_{n-1}(J)$.
Similarly, $q_n\in\CH^n_n(J)$.
It is also well known that
$p_1=-\th$ (see e.g. \cite{Bmain}, Prop. 2.3, or \cite{P-av}, 17.2 and 17.3). 

Recall that we have defined the 
tautological subring $\TT\CH(J)_{\Q}\sub\CH(J)_{\Q}$ as the smallest subring containing $[C]$ and
closed under $S$ and under all the pull-back operations $[n]^*$. Equivalently, this
is the smallest subring closed under $S$ and containing all classes $[C]_n$.
This immediately implies that all classes $p_n$ and $q_n$ belong to $\TT\CH(J)_{\Q}$.

Let us consider the element 
$$\eta:=\iota_*K/2+[0]\in\CH^g(J)_{\Q}.$$ 
From the Riemann's Theorem we get
\begin{equation}\label{thetaCeq}
\eta=\th\cdot [C]
\end{equation}
(see e.g. \cite{P-av}, Thm. 17.4).
Hence, $\th\cdot [C]_n=\eta_n$ and we have
$$q_n=S(\eta_n).$$
Note that for every point $x\in J$ we have
$$S([x])=\exp(c_1(\LL_x))=\exp(\th_x-\th), $$
where $\LL_x=\LL|_{J\times x}$ and $\th_x=[\Th+x]$.
Hence, we can rewrite the definition of $q_n$ for $n>0$ as follows:
$$q_n=\frac{1}{2}\sum_{i=1}^{2g-2}c_1(\LL_{x_i})^n,$$
where $(x_i)$ are points on $C$ such that $K=x_1+\ldots+x_{2g-2}$.
In particular, 
$$q_1=\frac{1}{2}c_1(\LL_{\kappa}),$$ 
where $\kappa\in J$ is the point corresponding to
$K(-(2g-2)x_0)$. 


\noindent
{\it Proof of Corollary \ref{zero-cor}.}
Theorem \ref{mainthm} implies that $\TT\CH(J)_{\Q}$ is generated with respect to the Pontryagin
product by the classes $([C]_n)$ and $(\eta_n)$. Therefore, the group of tautological $0$-cycles
is generated with respect to this product by the 
classes $(\eta_n)$, or equivalently, by the classes $([n]_*\eta)$.
\ed

The action of $\ssl_2$ on $\CH(J)_{\Q}$ (in fact, on the motive
of $J$) mentioned in the introduction is generated by the operators
$$e(a)=p_1\cdot a=-\th\cdot a,$$
$$f(a)=-[C]_0*a,$$
$$h(a)=(2n-s-g)a\text{ for }a\in\CH^n_s$$
(the operators $e$ and $f$ differ from those of \cite{K} by the sign).
In fact, this action is induced by an algebraic action of the group $\SL_2$, so that
the Fourier transform corresponds to the action of the matrix
$\left(\matrix 0 & 1 \\ -1 & 0\endmatrix\right)$ (see \cite{P-Wr}, Thm.~5.1).
This leads to the formula 
$$S=\exp(e)\exp(-f)\exp(e)$$
that can also be checked directly (see \cite{Bmain}, (1.7)).

\section{Commutation relations}\label{com-rel-sec}

Let us consider the following family of binary operations on $\CH(J)$:
$$a*_n b=(p_1+p_2)_*(c_1(\LL)^n\cdot p_1^*a\cdot p_2^*a), \ n\ge 0,$$
where $a,b\in\CH(J)$. Note that $a*_0 b=a*b$ is the usual Pontryagin
product.

\begin{lem}\label{adthetalem} One has 
$$a*_{n+1}b=(\th\cdot a)*_n b+a*_n(\th\cdot b)-\th\cdot(a*_n b).$$
\end{lem}

\Pf . This follows immediately from the identity
$c_1(\LL)=p_1^*\th+p_2^*\th-(p_1+p_2)^*\th$.
\ed

\begin{lem} If $a\in\CH^{p_1}_{s_1}$, $b\in\CH^{p_2}_{s_2}$, then
$a*_n b\in\CH^{p_1+p_2+n-g}_{s_1+s_2}$.
\end{lem}

\Pf . Since $\th\cdot\CH^p_s\sub\CH^{p+1}_s$, the assertion follows from Lemma
\ref{adthetalem} by induction in $n$.
\ed

For every $a\in\CH(J)$ and $n\ge 0$ let us denote by $A_n(a)$ the operator
$b\mapsto a*_n b$ on $\CH(J)$. For $n<0$ we set $A_n(a)=0$.
Note that Lemma \ref{adthetalem} is equivalent to the following identity
\begin{equation}\label{adeeq}
[e,A_n(a)]=A_{n+1}(a)-A_n(\th\cdot a),
\end{equation}
where $e$ is the operator of the $\ssl_2$-action (see section \ref{prelim-sec}).

\begin{lem}\label{vanishlem} For every $s\ge 0$ one has
$A_n(\eta_s)=0$ for $n>s$ and
$A_n([C]_s)=0$ for $n>s+2$.
\end{lem}

\Pf . 
We start by observing that the operator $f=-A_0([C]_0)$ commutes with $A_0([C]_s)$ and
with $A_0(\eta_s)$.
Also, 
$$[h,A_0([C]_s)]=(-s-2)A_0([C]_s),\ [h,A_0(\eta_s)]=-sA_0(\eta_s).$$
Hence, $A_0([C]_s)$ and $A_0(\eta_s)$ are lowest weight vectors with respect to the adjoint
action of $\ssl_2$ on $\End(\CH(J))$ of weights $-(s+2)$ and $-s$, respectively. It follows
that 
$$\ad(e)^n(A_0([C]_s))=0$$
for $n>s+2$ and
$$\ad(e)^n(A_0(\eta_s))=0$$
for $n>s$. Using equalities (\ref{adeeq}) and
(\ref{thetaCeq}) we find by induction in $n$ that 
$$\ad(e)^n(A_0(\eta_s))=A_n(\eta_s),$$
$$\ad(e)^n(A_0([C]_s))=A_n([C]_s)-nA_{n-1}(\eta_s).$$
The first equality implies that $A_n(\eta_s)=0$ for $n>s$.
Together with the second equality this implies that $A_n([C]_s)=0$ for
$n>s+2$.
\ed

\begin{lem}\label{Fourlem} One has
$$A_s(\eta_s)(x)=s!\cdot q_s\cdot x,$$
$$A_{s+2}([C]_s)(x)=(s+2)!\cdot p_{s+1}\cdot x.$$
\end{lem}

\Pf . 
As we have seen in the previous proof, the operator $A_0([C]_s)$ (resp., $A_0(\eta_s)$)
is a lowest weight vector of weight $-(s+2)$ (resp., $-s$) with respect to the $\ssl_2$-action.
Since the Fourier transform $S$ is given by the action of
$\left(\matrix 0 & 1 \\ -1 & 0\endmatrix\right)\in\SL_2$, it follows that 
$$\ad(e)^{s+2}(A_0([C]_s))=\la_s\cdot S A_0([C]_s) S^{-1},$$
$$\ad(e)^{s}(A_0(\eta_s))=\mu_s\cdot S A_0(\eta_s) S^{-1}$$
for some nonzero constants $\la_s$, $\mu_s$. 
But $\ad(e)^s(A_0(\eta_s))=A_s(\eta_s)$ and
$$\ad(e)^{s+2}(A_0([C]_s))=A_{s+2}([C]_s)-(s+2)A_{s+1}(\eta_s)=A_{s+2}([C]_s)$$
as we have seen in the proof of Lemma \ref{vanishlem}.
Hence,
$$A_s(\eta_s)(a)=\mu_s\cdot S(\eta_s)\cdot a,$$
$$A_{s+2}([C]_s)(a)=\la_s\cdot S([C]_s)\cdot a.$$
Setting $a=1$ we get
$$(p_1+p_2)_*(c_1(\LL)^s\cdot p_1^*\eta_s)=\mu_s\cdot S(\eta_s),$$
$$(p_1+p_2)_*(c_1(\LL)^{s+2}\cdot p_1^*[C]_s)=\la_s\cdot S([C]_s).$$
Making the change of variables $(x,y)\mapsto (x,x+y)$ in the formula
defining the Fourier transform and using the theorem
of the cube we get 
\begin{equation}\label{Fourformula}
S(a)=(p_1+p_2)_*(\exp(c_1(\LL)+p_1^*\De^*c_1(\LL))\cdot p_1^*a)=
(p_1+p_2)_*(\exp(c_1(\LL)+2p_1^*\th)\cdot p_1^*a),
\end{equation}
where $a\in\CH(J)_{\Q}$.
Applying this to $a=\eta_s$ and keeping in mind that $\th\cdot\eta_s=0$ we get
$$S(\eta_s)=(p_1+p_2)_*(\exp(c_1(\LL))\cdot p_1^*\eta_s).$$
Since $S(\eta_s)\in\CH^{g-s}$, this implies that
$$S(\eta_s)=(p_1+p_2)_*(\frac{c_1(\LL)^s}{s!}\cdot p_1^*\eta_s),$$
so $\mu_s=s!$.
Similarly, applying (\ref{Fourformula}) to $a=[C]_s$ and using the fact
that $\th^2\cdot[C]_s=0$ and $S([C]_s)\in\CH^{g-s-1}$ we obtain
$$S([C]_s)=(p_1+p_2)_*(\frac{c_1(\LL)^{s+2}}{(s+2)!}\cdot p_1^*[C]_s)+
2 (p_1+p_2)_*(\frac{c_1(\LL)^{s+1}}{(s+1)!}\cdot p_1^*(\th\cdot [C]_s)).$$
It remains to observe that the second term is proportional to
$A_{s+1}(\eta_s)(1)$, hence it vanishes by Lemma \ref{vanishlem}.
\ed

\begin{lem}\label{commutatorlem} For $n_1,n_2\ge 0$ and $a_1,a_2,b\in\CH(J)_{\Q}$ one has 
\begin{align*}
&[A_{n_1}(a_1),A_{n_2}(a_2)](b)=\\
&\sum_{i\ge 1}
(p_1+p_2+p_3)_*[\left({n_1\choose i}p_{13}^*c_1(\LL)^{n_1-i}p_{23}^*c_1(\LL)^{n_2}-
{n_2\choose i}p_{13}^*c_1(\LL)^{n_1}p_{23}^*c_1(\LL)^{n_2-i}\right)\\
&\cdot p_{12}^*c_1(\LL)^i\cdot p_1^*a_1\cdot p_2^*a_2\cdot p_3^*b],
\end{align*}
where $p_{ij}:J\times J\times J\to J\times J$ and $p_i:J\times J\times J\to J$ are the natural
projections.
\end{lem}

\Pf . Using the projection formula we find
$$a_1*_{n_1}(a_2*_{n_2}b)=(p_1+p_2+p_3)_*\left(
(p_1,p_2+p_3)^*c_1(\LL)^{n_1}\cdot p_{23}^*c_1(\LL)^{n_2}\cdot 
p_1^*a_1\cdot p_2^*a_2\cdot p_3^*b\right).$$
Similarly,
\begin{align*}
&a_2*_{n_2}(a_1*_{n_1}b)=(p_1+p_2+p_3)_*\left(
(p_1,p_2+p_3)^*c_1(\LL)^{n_2}\cdot p_{23}^*c_1(\LL)^{n_1}\cdot 
p_1^*a_2\cdot p_2^*a_1\cdot p_3^*b\right)=\\
&(p_1+p_2+p_3)_*\left(
(p_2,p_1+p_3)^*c_1(\LL)^{n_2}\cdot p_{13}^*c_1(\LL)^{n_1}\cdot 
p_1^*a_1\cdot p_2^*a_2\cdot p_3^*b\right).
\end{align*}
It remains to use the equalities 
$$(p_1,p_2+p_3)^*c_1(\LL)=p_{12}^*c_1(\LL)+p_{13}^*c_1(\LL),\
(p_2,p_1+p_3)^*c_1(\LL)=p_{12}^*c_1(\LL)+p_{23}^*c_1(\LL).
$$
\ed

Note that from the above lemma (or directly from the definition) one can immediately see
that for $a_1,a_2\in\CH^g(J)_{\Q}$ the operators $A_{n_1}(a_1)$ and $A_{n_2}(a_2)$
commute. Hence, $[A_{n_1}(\eta_{s_1}),A_{n_2}(\eta_{s_2})]=0$.

\begin{thm} One has the following commutation relations
\begin{align*}
&[A_{n_1}([C]_{s_1}),A_{n_2}([C]_{s_2})]=\\
&\left(n_1\cdot {s_1+s_2-n_1-n_2+3\choose s_1-n_1+2}
-n_2\cdot {s_1+s_2-n_1-n_2+3\choose s_2-n_2+2}\right)
A_{n_1+n_2-1}([C]_{s_1+s_2})-\\
&2\cdot\left({n_1\choose 2}{s_1+s_2-n_1-n_2+2\choose s_1-n_1+2}
-{n_2\choose 2}{s_1+s_2-n_1-n_2+2\choose s_2-n_2+2}\right)A_{n_1+n_2-2}(\eta_{s_1+s_2}),
\end{align*}
\begin{align*}
&[A_{n_1}([C]_{s_1},A_{n_2}(\eta_{s_2})]=\\
&\left(n_1{s_1+s_2-n_1-n_2+1\choose s_1-n_1+2}-n_2{s_1+s_2-n_1-n_2+1\choose s_2-n_2}\right)
A_{n_1+n_2-1}(\eta_{s_1+s_2}).
\end{align*}
\end{thm}

\Pf . Since $[m]_*[C]=\sum_{s\ge 0} m^{s+2}[C]_s$ (resp.,
$[m]_*\eta_s=\sum_{s\ge 0}m^s \eta_s$), the first (resp., the second) commutator is
the coefficient with $m_1^{s_1+2}m_2^{s_2+2}$ (resp., $m_1^{s_1+2}m_2^{s_2}$)
in the commutator
$$[A_{n_1}([m_1]_*[C]),A_{n_2}([m_2]_*[C])]\text{      (resp., }[A_{n_1}([m_1]_*[C]),A_{n_2}([m_2]_*\eta)]).$$ 
Let us apply Lemma \ref{commutatorlem} to compute these commutators.
Taking into account the formula for $\LL|_{C\times C}$ we find
\begin{align*}
&[A_{n_1}([m_1]_*[C]),A_{n_2}([m_2]_*[C])](b)=\\
&\sum_{i\ge 1}
(p_1+p_2+p_3)_*[\left({n_1\choose i}p_{13}^*c_1(\LL)^{n_1-i}p_{23}^*c_1(\LL)^{n_2}-
{n_2\choose i}p_{13}^*c_1(\LL)^{n_1}p_{23}^*c_1(\LL)^{n_2-i}\right)\\
&\cdot m_1^im_2^i
p_{12}^*\left([m_1]\times [m_2])_*(\iota\times \iota)_*([\De_C]-[x_0\times C]-[C\times x_0])^i\right)
\cdot p_3^*b].
\end{align*}
Therefore, the sum has only two terms, $T_1$ and $T_2$, corresponding to $i=1$ and $i=2$.
We have
\begin{align*}
&T_1=
(p_1+p_2+p_3)_*([n_1p_{13}^*c_1(\LL)^{n_1-1}p_{23}^*c_1(\LL)^{n_2}-\\
&n_2p_{13}^*c_1(\LL)^{n_1}p_{23}^*c_1(\LL)^{n_2-1}]
m_1m_2
p_{12}^*\left((m_1,m_2)_*[C]-[0]\times [m_2]_*[C]-[m_1]_*[C]\times [0]\right)
\cdot p_3^*b).
\end{align*}
Note that since we are only interested in the coefficient with $m_1^{s_1+2}m_2^{s_2+2}$
we can discard the terms linear in $m_1$ or $m_2$. Therefore, we can replace $T_1$ with
$$T'_1=
((m_1+m_2)p_1+p_2)_*\left(m_1m_2[n_1m_1^{n_1-1}m_2^{n_2}-n_2m_1^{n_1}m_2^{n_2-1}]
c_1(\LL)^{n_1+n_2-1}\cdot p_1^*[C]\cdot p_2^*b\right).$$
Using the formula
$$c_1(\LL)=\frac{([m_1+m_2]\times\id)^*c_1(\LL)}{m_1+m_2}$$
we obtain
\begin{align*}
&T'_1=\frac{n_1m_1^{n_1}m_2^{n_2+1}-n_2m_1^{n_1+1}m_2^{n_2}}{(m_1+m_2)^{n_1+n_2-1}}
(p_1+p_2)_*\left(c_1(\LL)^{n_1+n_2-1}\cdot p_1^*[m_1+m_2]_*[C]\cdot p_2^*b\right)=\\
&(n_1m_1^{n_1}m_2^{n_2+1}-n_2m_1^{n_1+1}m_2^{n_2})\sum_s (m_1+m_2)^{s-n_1-n_2+3}
A_{n_1+n_2-1}([C]_s)(b).
\end{align*}
Let us observe that by Lemma \ref{vanishlem} we can restrict the summation to $s$
such that $n_1+n_2-1\le s+2$, i.e., $s\ge n_1+n_2-3$.

On the other hand, using the formula $[\De_C]^2=-\De_*K$ we obtain
$$(\iota\times \iota)_*([\De_C]-[x_0\times C]-[C\times x_0])^2=-2\De_*\eta.$$
Hence,
\begin{align*}
&T_2=-2
(p_1+p_2+p_3)_*[\left({n_1\choose 2}p_{13}^*c_1(\LL)^{n_1-2}p_{23}^*c_1(\LL)^{n_2}-
{n_2\choose 2}p_{13}^*c_1(\LL)^{n_1}p_{23}^*c_1(\LL)^{n_2-2}\right)\\
&\cdot m_1^2m_2^2
p_{12}^*((m_1,m_2)_*\iota_*\eta)\cdot p_3^*b].
\end{align*}
We can rewrite this as
\begin{align*}
&-T_2=2\frac{{n_1\choose 2}m_1^{n_1}m_2^{n_2+2}-{n_2\choose 2}m_1^{n_1+2}m_2^{n_2}}
{(m_1+m_2)^{n_1+n_2-2}}
(p_1+p_2)_*\left(c_1(\LL)^{n_1+n_2-2}\cdot p_1^*([m_1+m_2]_*\eta)\cdot p_2^*b\right)=\\
&2\left({n_1\choose 2}m_1^{n_1}m_2^{n_2+2}-{n_2\choose 2}m_1^{n_1+2}m_2^{n_2}\right)
\sum_s (m_1+m_2)^{s-n_1-n_2+2}A_{n_1+n_2-2}(\eta_s)(b)
\end{align*}
Again by Lemma \ref{vanishlem} we can restrict the summation to $s$ such that
$s\ge n_1+n_2-2$. Now the required formula for 
$[A_{n_1}([C]_{s_1}),A_{n_2}([C]_{s_2})]$
follows easily by considering the coefficients
with $m_1^{s_1+2}m_2^{s_2+2}$ in $T'_1$ and $T_2$.

Following similar steps we can write
\begin{align*}
&[A_{n_1}([m_1]_*[C]),A_{n_2}([m_2]_*\eta)](b)=\\
&(p_1+p_2+p_3)_*[\left(n_1p_{13}^*c_1(\LL)^{n_1-1}p_{23}^*c_1(\LL)^{n_2}-
n_2p_{13}^*c_1(\LL)^{n_1}p_{23}^*c_1(\LL)^{n_2-1}\right)\\
&\cdot m_1m_2
p_{12}^*\left((m_1,m_2)_*\eta-[0]\times [m_2]_*\eta\right)
\cdot p_3^*b].
\end{align*}
Since we are interested in the coefficient with $m_1^{s_1+2}m_2^{s_2}$, we
can discard the term linear in $m_1$. Hence, $[A_{n_1}([C]_{s_1}),A_{n_2}(\eta_{s_2})]$
is equal to the coefficient with $m_1^{s_1+2}m_2^{s_2}$ in 
\begin{align*}
&\frac{n_1m_1^{n_1}m_2^{n_2+1}-n_2m_1^{n_1+1}m_2^{n_2}}{(m_1+m_2)^{n_1+n_2-1}}
(p_1+p_2)_*\left(c_1(\LL)^{n_1+n_2-1}\cdot p_1^*[m_1+m_2]_*\eta\cdot p_2^*b\right)=\\
&(n_1m_1^{n_1}m_2^{n_2+1}-n_2m_1^{n_1+1}m_2^{n_2})\sum_s (m_1+m_2)^{s-n_1-n_2+1}
A_{n_1+n_2-1}(\eta_s)(b),
\end{align*}
where the summation can be restricted to $s\ge n_1+n_2-1$ by Lemma \ref{vanishlem}.
This immediately implies the result.
\ed

Setting 
$$\wt{X}_{k,n}=k!\cdot A_n([C]_{k+n-2}),$$
\begin{equation}\label{Y-eq}
Y_{k,n}=k!\cdot A_n(\eta_{k+n})
\end{equation}
for $n\ge 0$, $k\ge 0$
we see that these operators satisfy the commutation relations
$$[\wt{X}_{k,n},\wt{X}_{k',n'}]=(nk'-n'k)\wt{X}_{k+k'-1,n+n'-1}-4\cdot
\left({n\choose 2}{k'\choose 2}-{n'\choose 2}{k\choose 2}\right)Y_{k+k'-2,n+n'-2},$$
$$[\wt{X}_{k,n},Y_{k',n'}]=(nk'-n'k)Y_{k+k'-1,n+n'-1},$$
where we set $\wt{X}_{k,n}=Y_{k,n}=0$ for $k<0$
(note that this convention agrees with Lemma \ref{vanishlem}).

Also, by Lemma \ref{Fourlem} we have
\begin{equation}\label{E0neq}
\frac{1}{n!}Y_{0,n}(a)=q_n\cdot a,
\end{equation}
\begin{equation}\label{D0neq}
\frac{1}{n!}\wt{X}_{0,n}(a)=p_{n-1}\cdot a.
\end{equation}

\begin{lem} One has
$$\wt{X}_{2,0}=-2f,$$
$$\wt{X}_{1,1}=-h+g\cdot\id,$$
$$\wt{X}_{0,2}=2e.$$
\end{lem}

\Pf . The first equality holds by the definition of $f$.
The third equality follows from (\ref{D0neq}).
It remains to use the relation
$$\frac{1}{4}\cdot[\wt{X}_{0,2},\wt{X}_{2,0}]=\wt{X}_{1,1}-Y_{0,0}=\wt{X}_{1,1}-g\cdot\id.$$
\ed

\begin{lem}\label{adeflem} 
One has the following relations
$$\ad(e)^n(Y_{k,0})=\frac{k!}{(k-n)!}\cdot Y_{k-n,n},$$
$$\ad(f)^n(Y_{0,k})=\frac{k!}{(k-n)!}\cdot Y_{n,k-n},$$
$$\ad(e)^n(\wt{X}_{k,0})=\frac{k!}{(k-n)!}\cdot \left(\wt{X}_{k-n,n}-n(k-n)Y_{k-n-1,n-1}\right),$$
$$\ad(f)^n(\wt{X}_{0,k})=\frac{k!}{(k-n)!}\cdot \left(\wt{X}_{n,k-n}-n(k-n)Y_{n-1,k-n-1}\right).$$
\end{lem}

\Pf . We have
$$[e,Y_{k,n}]=\frac{1}{2}[\wt{X}_{0,2},Y_{k,n}]=kY_{k-1,n+1},$$
$$[e,\wt{X}_{k,n}]=\frac{1}{2}[\wt{X}_{0,2},\wt{X}_{k,n}]=k\wt{X}_{k-1,n+1}-k(k-1)Y_{k-2,n}.$$
From this one can easily deduce the first and the third formulas by induction in $n$.
The other two are proved in the same way since our relations are skew-symmetric
with respect to switching $\wt{X}_{k,n}$ with $\wt{X}_{n,k}$ and $Y_{k,n}$ with $Y_{n,k}$.
\ed

Finally, we set 
\begin{equation}\label{X-eq}
X_{k,n}=\wt{X}_{k,n}-kn Y_{k-1,n-1}=
k!\cdot A_n([C]_{k+n-2})-(k-1)!\cdot A_{n-1}(\eta_{k+n-2}).
\end{equation}

\noindent
{\it Proof of Theorem \ref{com-rel-thm}.}
An easy computation shows that
$X_{k,n}$ and $Y_{k,n}$ satisfy relations \eqref{main-com-rel} and
\eqref{main-com-rel2}. It remains to check \eqref{FourXmn} and
\eqref{FourYmn}. We have 
$$SY_{m,0}S^{-1}=m!SA_0(\eta_m)S^{-1}=m!q_m=A_m(\eta_m)=Y_{0,m},$$
$$S\wt{X}_{m,0}S^{-1}=m!SA_0([C]_{m-2})S^{-1}=m!p_{m-1}=
A_m([C]_{m-2}=\wt{X}_{0,m}.$$
Since $SeS^{-1}=-f$ we immediately derive (\ref{FourYmn}) and (\ref{FourXmn})
from Lemma \ref{adeflem}.
\ed

\section{Proof of Theorem \ref{mainthm}.}

\noindent
(i) Let us denote by $\TT'\CH(J)_{\Q}\sub\CH(J)_{\Q}$ the subring generated by the classes
$(p_n)$ and $(q_n)$.
Consider the operator $f=-SeS^{-1}=-\wt{X}_{2,0}/2=-X_{2,0}/2$ on $\CH(J)_{\Q}$. We are going to show
that it preserves $\TT'\CH(J)_{\Q}$.
Note that
$$[f,p_n]=-\frac{1}{(n+1)!}[X_{2,0}/2,X_{0,n+1}]=\frac{1}{n!}X_{1,n},$$
\begin{equation}\label{Upp-eq}
[[f,p_n],p_m]=\frac{1}{(m+1)!n!}[X_{1,n-1},X_{0,m+1}]=-\frac{1}{m!n!}X_{0,m+n}=
-{m+n\choose m}p_{m+n-1},
\end{equation}
\begin{equation}\label{Upq-eq}
[[f,p_n],q_m]=\frac{1}{n!m!}[X_{1,n},Y_{0,m}]=-\frac{1}{n!(m-1)!}Y_{0,m+n-1}=
-{m+n-1\choose m-1}q_{m+n-1},
\end{equation}
\begin{equation}\label{Uqq-eq}
\text{and    }[[f,q_n],q_m]=0.
\end{equation}
On the other hand, from the definition of $q_n$ we derive
$$q_n=-Se([C]_n)=fS([C]_n)=f(p_{n+1}).$$
Since $f(1)=0$, this gives
\begin{equation}\label{Up-eq}
[f,p_n](1)=f(p_n)=q_{n-1}.
\end{equation}
Also, 
$$0=-Se(\eta_n)=fS(\eta_n)=f(q_n),$$
so 
\begin{equation}\label{Uq-eq}
[f,q_n](1)=0.
\end{equation}
We claim that these relations imply that
for any polynomial $F$ in $(p_n)$ and $(q_n)$ one has
\begin{equation}\label{U-eq}
f(F)=-\DD(F),
\end{equation}
where $\DD$ is the differential operator defined in the formulation of the theorem.
Indeed, from relations \eqref{Upp-eq}-\eqref{Uqq-eq} we see that
$$[[f,x],y](F)=-[[\DD,x],y](F),$$
where $x$ and $y$ are any of the classes $(p_n)$ or $(q_n)$. 
By induction in the degree this implies that
$$[f,p_n](F)=-[\DD,p_n](F)\text{  and  } [f,q_n](F)=-[\DD,q_n](F),$$
where the base of induction follows from relation \eqref{Up-eq} and \eqref{Uq-eq}.
Finally, another induction in degree proves \eqref{U-eq}.

Thus, we proved the operator $f\in\ssl_2$ preserves $\TT'\CH(J)_{\Q}$ and acts on it by the
differential operator $-\DD$. Since $\TT'\CH(J)_{\Q}$ is also closed under the operator $e\in\ssl_2$,
it follows that $\TT'\CH(J)_{\Q}$ is closed under the Fourier transform 
$S=\exp(e)\exp(-f)\exp(e)$.
Therefore, $\TT'\CH(J)_{\Q}$ coincides with the tautological subring $\TT\CH(J)_{\Q}$. 

\noindent
(ii) For $m=0$ the operator $X_{0,n}$ (resp., $Y_{0,n}$) acts as multiplication by
$n!p_{n-1}$ (resp., $n!q_n$). In particular, they preserve $\TT\CH(J)_{\Q}$.
The general case follows from this by induction in $m$ using commutation relations
$$[X_{2,0},X_{m,n}]=-2n X_{m+1,n-1},\ \ [X_{2,0},Y_{m,n}]=-2n Y_{m+1,n-1}$$
together with the fact that $X_{2,0}/2$ acts on $\TT\CH(J)_{\Q}$ via the operator $\DD$.
\ed

\section{Some relations in $\TT\CH(J)_{\Q}$}\label{ex-sec}

Let us denote by $\Q[q]\sub\TT\CH(J)_{\Q}$ the subring generated
by the classes $(q_n)$.
First, we collect some general observations in the following

\begin{prop}\label{taut-g-prop}
(i) $\oplus_s\TT\CH(J)^s_s=\Q[q]$.
 
\noindent
(ii) $\TT\CH(J)^g_s=p_1^{g-s}\cdot \TT\CH(J)^s_s=p_1^{g-s}\cdot(\CH^s\cap\Q[q])$.

\noindent
(iii) $\TT\CH(J)^g_g=0$. 
\end{prop}

\Pf .(i) This follows from the fact that $q_n\in\CH(J)^n_n$ and $p_n\in\CH(J)^n_{n-1}$.

\noindent
(ii) Since $f$ acts on $\TT\CH(J)_{\Q}$ by the operator $-\DD$, it preserves
the subring $\Q[p_1,q]$ generated by $p_1$ and $(q_n)$. Hence,
the Fourier transform $S$ also preserves this subring.
But $\TT\CH(J)^g_s=S(\CH^s_s)$, so the assertion follows from (i). 

\noindent
(iii) It is enough to prove that $q_{n_1}\ldots q_{n_k}=0$ for $n_1+\ldots+n_k=g$.
We can use induction in $k$. The base of induction follows from
$$q_g=-\DD(p_{g+1})=0.$$
Assume the assertion holds for $k-1$. Then for $n_1+\ldots+n_k=g$ we have
$$0=-\DD(q_{n_1}\ldots q_{n_{k-1}}p_{n_k+1})=q_{n_1}\ldots q_{n_{k-1}}q_{n_k},$$
since all the other terms vanish by the induction assumption.
\ed

Part (ii) of the above proposition implies that for every $n_1+\ldots+n_k+m_1+\ldots+m_l=g$ we have
a relation of the form
$$p_{n_1}\ldots p_{n_k}q_{m_1}\ldots q_{m_l}=p_1^{k}f(q).$$
The simplest example of such a relation is
\begin{equation}\label{pg-eq}
p_g=p_1q_{g-1}
\end{equation}
obtained by applying $\DD$ to $p_1p_g=0$.
Similarly, applying $\DD$ to $p_1p_iq_{g-i}=0$ we get
\begin{equation}\label{piq-eq}
p_iq_{g-i}=p_1q_{i-1}q_{g-i}-{g-1\choose i}p_1q_{g-1}.
\end{equation}

\begin{prop}
The ring $\TT\CH(J)_{\Q}$ is generated over $\Q$ 
by the classes $(p_n)_{n<g/2+1}$ and $(q_n)_{n<\frac{g+1}{2}}$. Furthermore,
for $n\ge g/2+1$ the class $p_n$ belongs to the ideal generated by $(q_i)_{i\ge 1}$.
\end{prop}

\Pf . First, let us prove that for $n\ge\frac{g+1}{2}$ the class $q_n$ can be expressed in
terms of $(q_i)_{i<n}$. The idea is to represent $n$ in the form $n=n_1+\ldots+n_k-k$,
where $n_i\ge 2$ for all $i$ and $n_1+\ldots+n_k\ge g+1$. Then it is enough to check
$0=\DD^k(p_{n_1}\ldots p_{n_k})$ is a polynomial of degree $n$ in $(q_i)_{i\le n}$ 
(where $\deg q_i=i$) that has a nonzero coefficient with $q_n$.
Note that $\DD$ preserves the subring generated by all the classes $(p_i)$, where $i\ge 2$,
together with all the classes $(q_i)$ and acts on this subring as the operator
$\DD'=\DD'_0-\DD'_1$, where
$$\DD'_0=\frac{1}{2}\sum_{m,n\ge 2}{m+n\choose n}p_{m+n-1}\del_{p_m}\del_{p_n}+
\sum_{m\ge 1,n\ge 2}{m+n-1\choose n}q_{m+n-1}\del_{q_m}\del_{p_n},$$
$$\DD'_1=\sum_{n\ge 2}q_{n-1}\del_{p_n}.$$
Let us consider two more gradings $\deg_p$ and $\deg_q$ on the algebra
of polynomials in $(p_i)$ and $(q_i)$ such that $\deg_p(p_i)=\deg_q(q_i)=1$ and
$\deg_p(q_i)=\deg_q(p_i)=0$. Since $\DD$ decreases $\deg_p$ by $1$, we obtain
that $\DD^k(p_{n_1}\ldots p_{n_k})$ is a polynomial in $(q_i)$. Furthermore,
since $\DD'_0$ preserves $\deg_q$ and $\DD'_1$ raises it by $1$, we have
$$\DD^k(p_{n_1}\ldots p_{n_k})=(\DD')^k(p_{n_1}\ldots p_{n_k})=
(\DD'_0)^k(p_{n_1}\ldots p_{n_k})-\la\cdot q_n+f(q_1,\ldots,q_{n-1}),$$
where 
$$\la\cdot q_n=\sum_{i=1}^k (\DD'_0)^{i-1}\DD'_1(\DD'_0)^{k-i}(p_{n_1}\ldots p_{n_k}).$$
Now it is clear from the explicit form of $\DD'_0$ and $\DD'_1$ that
$(\DD'_0)^k(p_{n_1}\ldots p_{n_k})=0$ and that $\la>0$.

Next, let us show that for $n\ge g/2+1$ the class $p_n$ belongs to the ideal generated by
$(q_i)_{i\ge 1}$ in the subring generated by $(p_i)_{i<n}$ over $\Q[q]$.
To this end we represent $n$ in the form $n=n_1+\ldots+n_k-k+1$, where $n_i\ge 2$ for all $i$
and $n_1+\ldots+n_k\ge g$. Now let us consider the class
$$a=\DD^{k-1}(p_{n_1}\ldots p_{n_k})=(\DD')^{k-1}(p_{n_1}\ldots p_{n_k}).$$
Note that if $n_1+\ldots+n_k=g$ then we have $p_{n_1}\ldots p_{n_k}\in\CH^g_{g-k}(J)$.
Hence, by Proposition \ref{taut-g-prop}(ii) the class 
$p_{n_1}\ldots p_{n_k}$ belongs to the ideal generated by $(q_i)$
in the subring $\Q[p_1,q]\sub\TT\CH(J)_{\Q}$. Therefore the same is true about the class $a$ 
(and if $n_1+\ldots+n_k>g$ then $a=0$).
On the other hand, we can write 
$\DD'=\DD'_p+\DD'_q$, where
$$\DD'_p=\frac{1}{2}\sum_{m,n\ge 2}{m+n\choose n}p_{m+n-1}\del_{p_m}\del_{p_n},$$
$$\DD'_q=\sum_{m\ge 1,n\ge 2}{m+n-1\choose n}q_{m+n-1}\del_{q_m}\del_{p_n}
-\sum_{n\ge 2}q_{n-1}\del_{p_n}.$$
Since $\deg_p(a)=1$ and since the image of $\DD'_q$ is contained in the ideal generated by $(q_i)$,
 we obtain
$$a=(\DD')^{k-1}(p_{n_1}\ldots p_{n_k})=\mu\cdot p_n+a',$$
where
$$\mu\cdot p_n=(\DD'_p)^{k-1}(p_{n_1}\ldots p_{n_k})$$
and $a'$ is a linear combination of $f_i(q)p_{n-i}$ with $i\ge 1$ and $f_i(q)\in\Q[q]$.
It is clear from the formula for $\DD'_p$ that $\mu>0$, so our assertion follows.
\ed

For example, if $g=2$ then $\TT\CH(J)_{\Q}$ is generated by $p_1$ and $q_1$.
In fact, in this case we have $q_2=q_1^2=0$ (by Proposition \ref{taut-g-prop})
and $p_2=p_1q_1$ (by \eqref{pg-eq}).
For $g=3$ the above proposition states that $\TT\CH(J)_{\Q}$ is generated by $p_1$, $p_2$
and $q_1$. Indeed, first we see that $q_3=q_1q_2=q_1^3=0$ and $p_3=p_1q_2$. Also, applying
$\DD^2$ to $p_2^2=0$ we derive the relation $q_2=q_1^2/4$ (and hence $p_3=p_1q_1^2/4$).
In addition, \eqref{piq-eq} gives $p_2q_1=\frac{3}{4} p_1q_1^2$.

Finally, let us show that the classes of the subvarieties of special divisors $W_d\sub J$ can
be expressed as universal polynomials in $(p_n-q_n)$. Recall that $W_d$ is the image
of the natural map $C^{[d]}\to J: D\to \OO_C(D-dx_0)$. Let us set $w_i=[W_{g-i}]\in\CH^i(J)$.

\begin{prop} One has
$-p_k+q_k=N^k(w)$, where $N^k(w)$ is the $k$-th Newton polynomial in the classes $(w_i)$:
$N^k(w)=\frac{1}{k!}\sum_{i=1}^g\la_i^g$, where $(\la_i)$ are roots of 
$\la^g-w_1\la^{g-1}+\ldots+(-1)^gw_g=0$.
\end{prop}

\Pf . We will use the formula 
\begin{equation}\label{wi-for}
(-1)^iw_i=c_i(\de^*E_d),
\end{equation}
where $E_d$ is the vector bundle defined for $d\gg 0$ by $E_d=S(\OO_C(dx_0))$ and
$\de:J\to J$ is the involution $x\mapsto\kappa-x$, where $\kappa\in J$ corresponds to
the line bundle $K(-(2g-2)x_0)$.
(see \cite{P-av}, sec.~19.5).
On the other hand, using Grothendieck-Riemann-Roch formula one can easily find that
$$\ch_i(E_d)=p_i-q_i$$
for $i>0$ (where $\ch_i$ is the component of the Chern character of codimension $i$).
Hence,
$$\ch_i(\de^*E_d)=\de_*(p_i-q_i)=[\kappa]*[-1]_*(p_i-q_i)=(-1)^{i-1}[\kappa]*(p_i+q_i).$$
Next, we observe that
$$S([\kappa])=\exp(c_1(\LL_{\kappa}))=\exp(2q_1).$$
Using Theorem \ref{com-rel-thm} we obtain that
$$[\kappa]*a=\exp(2Y_{1,0})(a)$$
for $a\in\CH^*(J)_{\Q}$.
Now the formula for $Y_{1,0}$ from Theorem \ref{mainthm}(ii) implies that
$$[\kappa]*q_i=q_i,\ \ \ [\kappa]*p_i=p_i-2q_i.$$
Therefore,
$$\ch_i(\de^*E_d)=(-1)^{i-1}[\kappa]*(p_i+q_i)=(-1)^{i-1}(p_i-q_i).$$ 
Combining this with \eqref{wi-for} we immediately derive our assertion.
\ed

Thus, we have $p_1-q_1=-w_1$, $p_2-q_2=w_2-w_1^2/2$, $p_3-q_3=w_3/2+w_2w_1/2-w_1^3/6$,
etc. Hence, we can express the classes $(w_i)$ as polynomials in 
$(p_i-q_i)$.

\end{document}